\documentclass[11pt]{article}
\usepackage[francais]{babel}
\usepackage[latin1]{inputenc}
\usepackage[T1]{fontenc}
\usepackage{a4,amsmath,amsfonts}
\title{Géométrie, points entiers et courbes entières}
\author{Pascal Autissier}
\begin{document}

\maketitle

\newcommand{\D}{\displaystyle}

{\bf Abstract:} Let $X$ be a projective variety over a number field $K$ (resp.
over $\mathbb{C}$). Let $H$ be the sum of ``sufficiently many positive
divisors'' on $X$. We show that any set of quasi-integral points (resp. any
integral curve) in $X-H$ is not Zariski dense.\\


\section{Introduction}

Soient $K$ un corps de nombres et $S$ un ensemble fini de places de $K$. On
note $O_{K;S}$ l'anneau des $S$-entiers de $K$. On s'intéresse dans cet article
aux solutions dans $O_{K;S}^r$ de systèmes d'équations du type
$$\forall i\in\{1;\cdots;n\}\quad F_i(x_1;\cdots;x_r)=0\quad,\qquad(*)$$
où les $F_i$ sont des polynômes à $r$ variables et à coefficients dans
$O_{K;S}$. Pour formaliser cette étude, on utilise le langage de la géométrie
algébrique:\\

Désignons par $Y$ la ``variété algébrique'' sur $K$ définie par
$F_1=0,\cdots,F_n=0$. Tout ensemble de solutions de $(*)$ dans $O_{K;S}^r$
définit alors un ensemble (de points) $S$-entier sur $Y$.

Le problème est de donner des conditions géométriques suffisantes sur $Y$ pour
que tout ensemble $S$-entier soit non Zariski-dense dans $Y$.\\

Dans la suite, on se donne $Y$ sous la forme $Y=X-D$, où $X$ est une variété
projective sur $K$ de dimension $d\ge1$ et $D$ un diviseur effectif sur $X$.
L'esprit de la conjecture de Lang et Vojta ({\it cf} conjecture 4.2 de
\cite{Lang} p. 223) est qu'une telle condition suffisante s'exprime en
termes de ``positivité'' de $D$:\\

{\bf Conjecture (Lang, Vojta):} Soit $X$ une variété projective lisse sur $K$
de diviseur canonique ${\cal K}_X$. Soit $D$ un diviseur effectif sur $X$, à
croisements normaux. Posons $Y=X-D$. On suppose ${\cal K}_X+D$ gros (par
exemple ample) sur $X$. Alors tout ensemble $S$-entier sur $Y$ est
non Zariski-dense dans $Y$.\\

Les théorèmes de Siegel et de Faltings \cite{Fal1} montrent cette conjecture
lorsque $X$ est une courbe. Plus généralement, cet énoncé est connu de Faltings
\cite{Fal2} lorsque $X$ est une sous-variété de variété abélienne. Notons
cependant que la conjecture est encore largement ouverte: le cas où
$X=\mathbb{P}^2_K$ n'est par exemple pas connu.\\

Dans cet article, on démontre des cas particuliers de cette conjecture, lorsque
$D$ a ``suffisamment'' de composantes irréductibles. Plus précisément, on
prouve le résultat suivant ({\it cf} section 2 pour les définitions):\\

{\bf Théorème 1.1:} {\it Soit $X$ une variété projective sur $K$ de dimension
$d\ge2$. Soient $D_1;\cdots;D_{d\delta}$ des diviseurs effectifs presque amples
sur $X$ qui se coupent proprement deux à deux (avec $\delta\ge2$). On suppose
que toute intersection de $\delta+1$ quelconques d'entre eux est vide. Posons
$Y=X-D_1\cup\cdots\cup D_{d\delta}$. Alors $Y$ est arithmétiquement
quasi-hyperbolique. En particulier, tout ensemble ${\cal E}\subset Y(K)$
$S$-entier sur $Y$ est non Zariski-dense dans $Y$.}\\

Cet énoncé améliore un résultat récent de Levin ({\it cf} théorème 10.4A de
\cite{Lev}). En fait, Levin a besoin de
$\D2\Bigl\lfloor\frac{\delta+1}{2}\Bigr\rfloor d+1$ diviseurs au lieu de
$d\delta$.\\

On démontre en outre l'énoncé suivant (on prouve en fait un résultat plus
général):\\

{\bf Théorème 1.2:} {\it Soit $X$ une variété projective sur $K$ de dimension
$d\ge2$. Soient $D_1;\cdots;D_r$ des diviseurs effectifs amples sur $X$ qui se
coupent proprement (avec $r\ge2d$). Posons $\D L=\sum_{i=1}^rD_i$ et
$Y=X-D_1\cup\cdots\cup D_r$. On suppose que le diviseur $L-2dD_i$ est nef pour
tout $i\in\{1;\cdots;r\}$. Alors $Y$ est arithmétiquement
quasi-hyperbolique.}\\

L'hypothèse sur les $L-2dD_i$ est vérifiée lorsque les $D_i$ vivent dans un
cône ``suffisamment étroit'' du groupe de Néron-Severi de $X$. L'intérêt de ce
résultat réside dans le nombre (potentiellement linéaire en $d$) de diviseurs à
considérer.\\

Remarquons que les théorèmes 1.1 et 1.2 s'inscrivent bien dans le cadre de la
conjecture de Lang et Vojta, puisque si $X$ est lisse sur $K$ de diviseur
canonique ${\cal K}_X$ et les $D_i$ sont amples sur $X$, alors
${\cal K}_X+D_1+\cdots+D_r$ est ample sur $X$ dès que $r\ge d+2$ ({\it cf}
exemple 1.5.35 de \cite{Laza} p. 87).\\

Les démonstrations reposent sur une extension (théorème 3.3) de travaux de
Corvaja-Zannier \cite{CZ2} et de Levin \cite{Lev}, qui donne des conditions
géométriques de non-Zariski-densité des points $S$-entiers, et sur un bon choix
(théorème 4.4) de multiplicités associées aux diviseurs $D_i$.

L'ingrédient arithmétique principal est la version de Vojta \cite{Voj2} du
théorème du sous-espace de Schmidt \cite{Schm} et Schlickewei \cite{Schl}
(c'est un énoncé d'approximation diophantienne qui généralise le théorème de
Roth).\\

Par ailleurs, Vojta \cite{Voj1} a développé un ``dictionnaire'' entre la
géométrie diophantienne et la théorie de Nevanlinna: l'étude des points
$S$-entiers sur les variétés sur $K$ est mise en analogie avec l'étude des
courbes entières sur les variétés complexes.\\

Pour étayer ce dictionnaire, on montre aussi les énoncés qui ``correspondent''
aux théorèmes 1.1 et 1.2:\\

{\bf Théorème 1.3:} {\it Soit $X$ une variété complexe projective de dimension
$d\ge2$. Soient $D_1;\cdots;D_{d\delta}$ des diviseurs effectifs presque amples
sur $X$ qui se coupent proprement deux à deux (avec $\delta\ge2$). On suppose
que toute intersection de $\delta+1$ quelconques d'entre eux est vide. Posons
$Y=X-D_1\cup\cdots\cup D_{d\delta}$. Alors $Y$ est Brody quasi-hyperbolique. En
particulier, toute courbe entière $f:\mathbb{C}\rightarrow Y(\mathbb{C})$ est
d'image non Zariski-dense dans $Y$.}\\

{\bf Théorème 1.4:} {\it Soit $X$ une variété complexe projective de dimension
$d\ge2$. Soient $D_1;\cdots;D_r$ des diviseurs effectifs amples sur $X$ qui se
coupent proprement (avec $r\ge2d$). Posons $\D L=\sum_{i=1}^rD_i$ et
$Y=X-D_1\cup\cdots\cup D_r$. On suppose que le diviseur $L-2dD_i$ est nef pour
tout $i\in\{1;\cdots;r\}$. Alors $Y$ est Brody quasi-hyperbolique.}\\

Je remercie Antoine Chambert-Loir et Christophe Mourougane pour de fructueuses
discussions.\\

\section{Définitions et énoncés}

\subsection{Géométrie}

Soit $K$ un corps de caractéristique nulle.\\

{\bf Conventions:} On appelle variété sur $K$ tout schéma quasi-projectif et
géométriquement intègre sur $K$. Le mot ``diviseur'' sous-entend
``diviseur de Cartier''.\\

Soit $X$ une variété projective sur $K$ de dimension $d\ge1$. Lorsque $L$ est
un diviseur sur $X$ tel que $h^0(X;L)\ge1$, on désigne par ${\bf B}_L$ le lieu
de base de $\Gamma(X;L)$ et par
$\Phi_L:X-{\bf B}_L\rightarrow\mathbb{P}(\Gamma(X;L))$ le morphisme défini par
$\Gamma(X;L)$. Pour tout diviseur effectif $D$ sur $X$, on note $1_D$ la
section globale de ${\cal O}_X(D)$ qu'il définit.\\

{\bf Définition:} Un diviseur $L$ sur $X$ est dit {\bf libre} lorsque
${\bf B}_L$ est vide.\\

{\bf Définition:} Un diviseur $L$ sur $X$ est dit {\bf gros} lorsque
$\D\liminf_{n\rightarrow+\infty}\frac{1}{n^d}h^0(X;nL)>0$.\\

{\bf Définition:} Soit $L$ un diviseur gros sur $X$. On dit que $L$ est
{\bf presque ample} lorsqu'il existe un entier $n\ge1$ tel que $nL$ soit
libre.\\

{\bf Définition:} Un $\mathbb{R}$-diviseur $L$ sur $X$ est dit {\bf nef}
lorsque pour tout 1-cycle effectif $C$ sur $X$, on a $\bigl<L.C\bigr>\ge0$ (où
$\bigl<L.C\bigr>$ désigne le nombre d'intersection).\\

{\bf Définition:} Soient $D_1$ et $D_2$ deux diviseurs effectifs sur $X$. On
dit que $D_1$ et $D_2$ {\bf se coupent proprement} lorsque
${\cal O}_X(-D_1-D_2)={\cal O}_X(-D_1)\cap{\cal O}_X(-D_2)$.\\

{\bf Définition:} Plus généralement, soient $D_1;\cdots;D_r$ des diviseurs
effectifs sur $X$. On dit que $D_1;\cdots;D_r$ {\bf se coupent proprement}
lorsque pour toute partie $I$ non vide de $\{1;\cdots;r\}$, la section globale
$(1_{D_i})_{i\in I}$ de $\D\bigoplus_{i\in I}{\cal O}_X(D_i) $ est régulière
(autrement dit, pour tout $\D x\in\bigcap_{i\in I}D_i$, en notant $\varphi_i$
une équation locale de $D_i$ en $x$, les $(\varphi_i)_{i\in I}$ forment une
suite régulière de l'anneau local ${\cal O}_{X;x}$).\\

Soit $L$ un diviseur sur $X$ tel que $h^0(X;L)\ge1$. Soient $D_1;\cdots;D_r$
des diviseurs effectifs non nuls sur $X$ (avec $r\ge1$). Notons ${\cal P}$
l'ensemble des parties $I$ non vides de $\{1;\cdots;r\}$ telles que
$\D\bigcap_{i\in I}D_i$ soit non vide. Pour $I\in{\cal P}$,
$\underline{a}=(a_i)_i\in\mathbb{N}^I$ et $k\in\mathbb{N}^*$, on définit le
sous-espace vectoriel $V_{I;\underline{a};k}$ de $\Gamma(X;L)$ par
$$V_{I;\underline{a};k}=\sum_{\underline{b}}\Gamma\Bigl(X;L-\sum_{i\in I}b_iD_i\Bigr)\ \mbox{où la somme porte sur les}\ \underline{b}\in\mathbb{N}^I\ \mbox{tels que}\ \sum_{i\in I}a_ib_i\ge k\ .$$

{\bf Définition:} On pose $\D\nu(L;D_1;\cdots;D_r)=\inf_{I\in{\cal P}}\inf_{\underline{a}\in\mathbb{N}^I-\{0\}}\frac{\D\sum_{k\ge1}\dim V_{I;\underline{a};k}}{\D h^0(X;L)\sum_{i\in I}a_i}$.\\

On démontre à la section 4 le résultat suivant:\\

{\bf Théorème 2.1:} {\it On suppose $d\ge2$. Soient $D_1;\cdots;D_r$ des
diviseurs effectifs presque amples sur $X$ qui se coupent proprement deux à
deux; supposons que toute intersection de $\delta+1$ quelconques d'entre eux
est vide (avec $2\le\delta\le r$). Il existe alors
$(m_1;\cdots;m_r)\in\mathbb{N}^{*r}$ tel qu'en posant
$\D L=\sum_{i=1}^rm_iD_i$, on ait $\D\liminf_{n\rightarrow+\infty}\frac{1}{n}\nu(nL;m_1D_1;\cdots;m_rD_r)>\frac{r}{d\delta}$.}\\

Posons
$\D\lambda_d=\Bigl[1-\Bigl(1-\frac{1}{d}\Bigr)^{d+1}\Bigr]\frac{d}{d+1}$. On
prouve à la section 5.2 l'énoncé suivant:\\

{\bf Théorème 2.2:} {\it Soient $D_1;\cdots;D_r$ des diviseurs effectifs non
nuls et nefs sur $X$ qui se coupent proprement. On suppose que
$\D L=\sum_{i=1}^rD_i$ est ample. Soit $\theta>1$ un réel tel que le
$\mathbb{R}$-diviseur $L-d\theta D_i$ soit nef pour tout $i\in\{1;\cdots;r\}$.
On a alors la minoration
$$\liminf_{n\rightarrow+\infty}\frac{1}{n}\nu(nL;D_1;\cdots;D_r)\ge\lambda_d\theta\quad.$$}\\

{\bf Remarque:} On a $\D\lambda_d>\frac{1}{2}$ pour tout $d\ge2$. En outre,
$\lambda_d$ converge vers $1-e^{-1}$ lorsque $d$ tend vers $+\infty$.\\

\subsection{Hyperbolicité}

Lorsque $K$ est un corps de nombres et $S$ un ensemble fini de places de $K$,
on note $O_{K;S}$ l'anneau des $S$-entiers de $K$, {\it i.e.} l'ensemble des
$x\in K$ tels que $|x|_v\le1$ pour toute place finie $v\notin S$.\\

Soit $K$ un corps de nombres.\\

{\bf Définition:} Soient $Y$ une variété sur $K$, $K'$ une extension finie de
$K$ et $S$ un ensemble fini de places de $K'$. Un ensemble
${\cal E}\subset Y(K')$ est dit {\bf $S$-entier} sur $Y$ lorsqu'il existe un
$O_{K';S}$-schéma intègre et quasi-projectif ${\cal Y}$ de fibre générique
$Y_{K'}$ tel que ${\cal E}\subset{\cal Y}(O_{K';S})$.\\

{\bf Définition:} Soient $Y$ une variété sur $K$ et $K'$ une extension finie de
$K$. Un ensemble ${\cal E}\subset Y(K')$ est dit {\bf quasi-entier} sur $Y$
lorsqu'il existe un ensemble fini $S$ de places de $K'$ tel que ${\cal E}$ soit
$S$-entier sur $Y$.\\

{\bf Définition:} Soit $Y$ une variété sur $K$. On dit que $Y$ est
{\bf arithmétiquement quasi-hyperbolique} lorsqu'il existe un fermé $Z\neq Y$
tel que pour toute extension finie $K'$ de $K$ et tout ensemble quasi-entier
${\cal E}\subset Y(K')$ sur $Y$, l'ensemble ${\cal E}-Z(K')$ soit fini.\\

{\bf Définition:} Soit $Y$ une variété complexe. Une {\bf courbe entière} sur
$Y$ est une application holomorphe $f:\mathbb{C}\rightarrow Y(\mathbb{C})$ non
constante.\\

{\bf Définition:} Soit $Y$ une variété complexe. On dit que $Y$ est
{\bf Brody quasi-hyperbolique} lorsqu'il existe un fermé $Z\neq Y$ tel que
pour toute courbe entière $f$ sur $Y$, on ait
$f(\mathbb{C})\subset Z(\mathbb{C})$.\\

L'intérêt de la définition de $\nu$ réside dans les critères suivants:\\

{\bf Théorème (3.3):} {\it Soit $X$ une variété projective sur $K$ de dimension
$d\ge1$. Soient $D_1;\cdots;D_r$ des diviseurs effectifs non nuls sur $X$ qui
se coupent proprement deux à deux. Posons $Y=X-D_1\cup\cdots\cup D_r$. Soit
$n\ge1$ un entier. On suppose que le diviseur $\D L=n\sum_{i=1}^rD_i$ est libre
et gros sur $X$ et que $\nu(L;D_1;\cdots;D_r)>n$. Alors $Y$ est
arithmétiquement quasi-hyperbolique.}\\

{\bf Théorème (3.5):} {\it Soit $X$ une variété complexe projective de
dimension $d\ge1$. Soient $D_1;\cdots;D_r$ des diviseurs effectifs non nuls sur
$X$ qui se coupent proprement deux à deux. Posons $Y=X-D_1\cup\cdots\cup D_r$.
Soit $n\ge1$ un entier. On suppose que le diviseur $\D L=n\sum_{i=1}^rD_i$ est
libre et gros sur $X$ et que $\nu(L;D_1;\cdots;D_r)>n$. Alors $Y$ est Brody
quasi-hyperbolique.}\\

On en déduit le théorème 1.1 ---respectivement 1.3--- en appliquant le théorème
2.1 (avec $r=d\delta$) puis le théorème 3.3 ---respectivement 3.5---.

On en déduit de même les théorèmes 1.2 et 1.4 en appliquant le théorème 2.2.\\

\section{Démonstration des critères}

\subsection{Rappels}

Soit $X$ une variété complexe projective. Soit $L$ un faisceau inversible sur
$X$. On munit $L$ d'une métrique (continue) $\|\ \|$ et on pose
$\hat{L}=(L;\|\ \|)$.

Soit $f$ une courbe entière sur $X$. On définit la {\bf fonction
caractéristique} ${\rm T}_{\hat{L};f}:\mathbb{R}_+\rightarrow\mathbb{R}$ de $f$
relativement à $\hat{L}$ de la manière suivante:

On choisit une section rationnelle $s$ de $L$ définie et non nulle en $f(0)$.
Pour tout réel $r\ge0$, on pose
$${\rm T}_{\hat{L};f}(r)=\sum_{|z|\le r}\mu_z(f^*s)\ln\frac{r}{|z|}-\int_0^{2\pi}\ln\|s(f(re^{i\theta}))\|\frac{{\rm d}\theta}{2\pi}+\ln\|s(f(0))\|\quad,$$
où $\mu_z(f^*s)$ désigne l'ordre de $f^*s$ en $z\in\mathbb{C}$. Cela ne dépend
pas du choix de $s$.\\

Soient $K$ un corps de nombres et $X'$ une variété projective sur $K$. Soit
$L'$ un faisceau inversible sur $X'$. On munit $L'$ d'une métrique adélique
$(\|\ \|_v)_v$ et on pose $\hat{L'}=(L';(\|\ \|_v)_v)$ (pour des précisions sur
les métriques adéliques, on pourra consulter le paragraphe 1.2 de \cite{Zha}).

Soient $K'$ une extension finie de $K$ et $P\in X'(K')$. On définit la
{\bf hauteur} (normalisée) ${\rm h}_{\hat{L'}}(P)$ de $P$ relativement à
$\hat{L'}$ de la façon suivante:

On choisit une section rationnelle $s'$ de $L'$ définie et non nulle en $P$. On
pose
$${\rm h}_{\hat{L'}}(P)=-\frac{1}{[K':\mathbb{Q}]}\sum_v\ln\|s'(P)\|_v\quad,$$
où $v$ parcourt l'ensemble des places de $K'$. Ce réel ne dépend pas du choix
de $s'$.\\

\subsection{Cas arithmétique}

Commençons par un résultat facile d'algèbre linéaire:\\

{\bf Lemme 3.1:} {\it Soient $K$ un corps et $V$ un $K$-espace vectoriel de
dimension finie. Soit $({\cal F}_k)_{k\ge1}$ une suite décroissante de parties
de $V$ telle que ${\cal F}_k=\{0\}$ pour tout $k$ assez grand. Il existe alors
une base ${\cal B}$ de $V$ adaptée à la suite $({\cal F}_k)_{k\ge1}$,
{\it i.e.} ${\cal B}\cap{\cal F}_k$ est une base de ${\rm Vect}({\cal F}_k)$
pour tout $k\ge1$.}\\

{\it Démonstration:} Soit $m\ge1$ un entier tel que ${\cal F}_k=\{0\}$. On pose
${\cal B}_m=\emptyset$. On construit par récurrence une suite
$({\cal B}_m;\cdots;{\cal B}_1)$ de parties de $V$ de la manière suivante:

Pour $k\in\{1;\cdots;m-1\}$, on complète la partie libre ${\cal B}_{k+1}$ en
une base ${\cal B}_k$ de ${\rm Vect}({\cal F}_k)$ contenue dans ${\cal F}_k$.

Pour finir, on complète ${\cal B}_1$ en une base ${\cal B}$ de $V$. $\square$\\

Soient $K$ un corps de nombres et $X$ une variété projective sur $K$ de
dimension $d\ge1$. On va utiliser la version suivante du théorème du
sous-espace de Schmidt, Schlickewei et Vojta:\\

{\bf Proposition 3.2:} {\it Soit $L\in{\rm Pic}(X)$ libre et gros. Notons
$q=h^0(X;L)$. On munit $L$ d'une métrique adélique $(\|\ \|_v)_v$. Soient
$s_1;\cdots;s_N$ des sections non nulles engendrant $\Gamma(X;L)$. Soit
$\varepsilon>0$. Il existe alors un fermé $Z\neq X$ tel que pour toute
extension finie $K'$ de $K$ et tout ensemble fini $S$ de places de
$K'$, l'ensemble des points $P\in (X-Z)(K')$ vérifiant
$$\sum_{v\in S}\max_{J\in{\cal L}}\sum_{j\in J}\ln\|s_j(P)\|_v^{-1}\ge(q+q\varepsilon)[K':\mathbb{Q}]{\rm h}_{\hat{L}}(P)\qquad(1)$$
est fini, où ${\cal L}$ désigne l'ensemble des parties $J$ de $\{1;\cdots;n\}$
telles que $(s_j)_{j\in J}$ soit une base de $\Gamma(X;L)$.}\\

{\it Démonstration:} En posant $V=\Gamma(X;L)$, on a un morphisme
$\Phi_L:X\rightarrow\mathbb{P}(V)$ génériquement fini. Il existe donc un fermé
$Z_1\neq X$ tel que $\Phi_{L|X-Z_1}$ soit à fibres finies. On applique alors la
version de Vojta ({\it cf} théorème 0.3 et reformulation 3.4 de \cite{Voj2}) du
théorème du sous-espace:

Il existe une réunion finie $H$ de $K$-hyperplans de
$\mathbb{P}(V)\simeq\mathbb{P}^{q-1}_K$ telle que pour toute extension finie
$K'$ de $K$ et tout ensemble fini $S$ de places de $K'$, l'ensemble des points
$P\in (X-Z_1\cup\Phi_L^{-1}(H))(K')$ vérifiant $(1)$ est fini. $\square$\\

{\bf Remarque:} Vojta a en fait montré que l'on peut trouver un $Z$ indépendant
de $\varepsilon$, mais on n'en aura pas besoin dans la suite.\\

On montre ci-dessous une extension d'un résultat de Levin ({\it cf} théorème
8.3A de \cite{Lev}), lui-même inspiré de travaux de Corvaja et Zannier
({\it cf} théorème principal de \cite{CZ2} p. 707-708):\\

{\bf Théorème 3.3:} {\it Soient $D_1;\cdots;D_r$ des diviseurs effectifs non
nuls sur $X$ qui se coupent proprement deux à deux. Posons
$Y=X-D_1\cup\cdots\cup D_r$. Soit $m\ge1$ un entier. On suppose que le diviseur
$\D L=m\sum_{i=1}^rD_i$ est libre et gros sur $X$ et que
$\nu(L;D_1;\cdots;D_r)>m$. Alors $Y$ est arithmétiquement
quasi-hyperbolique.}\\

{\it Démonstration:} On procède en deux étapes: dans la première, on construit
un fermé $Z\neq X$ candidat à contenir ``presque tous les points entiers'';
dans la seconde, on prouve que $Y$ est arithmétiquement quasi-hyperbolique.\\

{\it \'Etape 1:} On pose $\D\varepsilon=\frac{1}{4m}(\nu(L;D_1;\cdots;D_r)-m)$
et $q=h^0(X;L)$, on choisit un entier $c\ge1$ tel que $h^0(X;L-cD_i)=0$ pour
tout $i\in\{1;\cdots;r\}$, et on fixe un entier
$\D b\ge\frac{cr}{m\varepsilon}$. Choisissons aussi une base ${\cal B}_0$ de
$\Gamma(X;L)$.\\

Désignons par ${\cal P}$ l'ensemble des parties $I$ non vides de
$\{1;\cdots;r\}$ telles que $\D\bigcap_{i\in I}D_i$ soit non vide. Soit
$I\in{\cal P}$. On note $\D\Delta_I$ l'ensemble des
$\underline{a}=(a_i)_i\in\mathbb{N}^I$ tels que $\D\sum_{i\in I}a_i=b$. Soit
$\underline{a}\in\Delta_I$. Pour $k\in\mathbb{N}^*$, on pose
$\D{\cal F}_k=\bigcup_{\underline{b}}\Gamma\Bigl(X;L-\sum_{i\in I}b_iD_i\Bigr)$ où la réunion porte sur les
$\underline{b}\in\mathbb{N}^I$ tels que $\D\sum_{i\in I}a_ib_i\ge k$, et on
note $V_{I;\underline{a};k}={\rm Vect}({\cal F}_k)$. Le lemme 3.1 fournit une
base ${\cal B}_{I;\underline{a}}$ de $\Gamma(X;L)$ adaptée à la suite
$({\cal F}_k)_{k\ge1}$.\\

On munit chaque faisceau ${\cal O}_X(D_i)$ d'une métrique adélique. Appliquons
le théorème du sous-espace (proposition 3.2) avec $\D\{s_1;\cdots;s_N\}={\cal B}_0\cup\bigcup_{I;\underline{a}}{\cal B}_{I;\underline{a}}$
(remarquons que cette réunion est finie puisque ${\cal P}$ et les $\Delta_I$ le
sont):

Il existe un fermé $Z\neq X$ tel que pour toute extension finie $K'$ de $K$ et
tout ensemble fini $S$ de places de $K'$, l'ensemble des points
$P\in (X-Z)(K')$ vérifiant l'inégalité $(1)$ est fini.\\

{\it \'Etape 2:} Soient $K'$ une extension finie de $K$ et $S$ un ensemble fini
de places de $K'$ contenant les places archimédiennes. Soit
${\cal E}\subset Y(K')$ un ensemble $S$-entier sur $Y$. Raisonnons par
l'absurde en supposant ${\cal E}-Z(K')$ infini. On choisit une suite injective
$(P_n)_{n\ge0}$ d'éléments de ${\cal E}-Z(K')$.

Quitte à extraire, on peut supposer (par compacité) que pour tout $v\in S$, la
suite $(P_{nv})_{n\ge0}$ converge dans $X(K'_v)$ vers un $y_v\in X(K'_v)$.\\

Pour tout $v\in S$, on note $I_v$ l'ensemble des $i\in\{1;\cdots;r\}$ tels que
$y_v\in D_i$. Quitte à extraire de nouveau, on peut supposer que pour tout
$v\in S$ tel que $I_v$ soit non vide et tout $i\in I_v$, la suite
$\D\biggl(\frac{\ln\|1_{D_i}(P_n)\|_v}{\D\sum_{j\in I_v}\ln\|1_{D_j}(P_n)\|_v}\biggr)_{n\ge0}$ converge vers un $t_{vi}\in[0;1]$. Remarquons que l'on a
$\D\sum_{i\in I_v} t_{vi}=1$.\\

\underline{Fait}: Soit $v\in S$. Il existe une base $(s_{1v};\cdots;s_{qv})$ de
$\Gamma(X;L)$ contenue dans $\{s_1;\cdots;s_N\}$ telle que l'on ait la
minoration suivante pour tout $n\ge0$:
$$-\sum_{k=1}^q\ln\|s_{kv}(P_n)\|_v\ge-(q+2q\varepsilon)\ln\|1_L(P_n)\|_v-O(1)\quad,\qquad(2)$$
où le $O(1)$ est indépendant de $n$.\\

Prouvons ce fait. Si $I_v$ est vide, on prend
$\{s_{1v};\cdots;s_{qv}\}={\cal B}_0$ et on obtient la minoration $(2)$ en
remarquant que $\ln\|1_L(P_n)\|_v=O(1)$ (puisque $y_v\notin L$).

On suppose maintenant $I_v$ non vide. On a donc $I_v\in{\cal P}$. Choisissons
un $\underline{a}_v=(a_{vi})_i\in\Delta_{I_v}$ tel que $|bt_{vi}-a_{vi}|\le1$
pour tout $i\in I_v$. On prend alors
$\{s_{1v};\cdots;s_{qv}\}={\cal B}_{I_v;\underline{a}_v}$. Vérifions que ce
choix convient.\\

Soit $s\in\Gamma(X;L)-\{0\}$. Pour tout
$i\in\{1;\cdots;r\}$, notons $\mu_i(s)$ le plus grand entier $\mu$ tel que le
diviseur ${\rm div}(s)-\mu D_i$ soit effectif. Puisque les diviseurs $D_i$ se
coupent proprement deux à deux, le diviseur
$\D{\rm div}(s)-\sum_{i\in I_v}\mu_i(s)D_i$ est encore effectif.
Ceci implique
$$-\ln\|s(P_n)\|_v\ge-\sum_{i\in I_v}\mu_i(s)\ln\|1_{D_i}(P_n)\|_v-O(1)\quad.$$

En remarquant que $\D t_{vi}>\frac{a_{vi}}{b}-\frac{2m\varepsilon}{cr}$ et que
$\mu_i(s)\le c$ pour tout $i\in I_v$, on a, par définition des $t_{vi}$:
$$-\mu_i(s)\ln\|1_{D_i}(P_n)\|_v\ge-\Bigl(\frac{a_{vi}}{b}\mu_i(s)-\frac{2m\varepsilon}{r}\Bigr)\sum_{j\in I_v}\ln\|1_{D_j}(P_n)\|_v$$
pour tout $n$ assez grand et tout $i\in I_v$.

On en déduit l'inégalité (pour tout $n\ge0$)
$$-\ln\|s(P_n)\|_v\ge-\Bigl(\frac{1}{b}\sum_{i\in I_v}a_{vi}\mu_i(s)-2m\varepsilon\Bigr)\sum_{j\in I_v}\ln\|1_{D_j}(P_n)\|_v-O(1)\quad.$$

On écrit cette inégalité pour $s=s_{kv}$, puis on somme sur $k$. En observant
que pour $i\in I_v$, on a
$$\begin{array}{rcl}
\D\sum_{k=1}^q\sum_{i\in I_v}a_{vi}\mu_i(s_{kv})&=&\D\sum_{\mu\ge1}\#\Bigl\{k\in\{1;\cdots;q\}\ \Big|\ \sum_{i\in I_v}a_{vi}\mu_i(s_{kv})\ge\mu\Bigr\}\\
&=&\D\sum_{\mu\ge1}\dim V_{I_v;\underline{a}_v;\mu}\ge\nu(L;D_1;\cdots;D_r)qb=(1+4\varepsilon)qbm\quad,\\
\end{array}$$
on trouve alors
$$-\sum_{k=1}^q\ln\|s_{kv}(P_n)\|_v\ge-(q+2q\varepsilon)m\sum_{j\in I_v}\ln\|1_{D_j}(P_n)\|_v-O(1)\quad.$$

Le fait énoncé $(2)$ s'en déduit en remarquant que $\ln\|1_{D_j}(P_n)\|_v=O(1)$
pour tout $j\notin I_v$.\\

Maintenant, l'ensemble ${\cal E}$ est $S$-entier sur $Y$, donc pour tout
$n\ge0$, on a
$$[K':\mathbb{Q}]{\rm h}_{\hat{L}}(P_n)=-\sum_{v\in S}\ln\|1_L(P_n)\|_v+O(1)\quad.$$

En utilisant la minoration $(2)$, on obtient (pour tout $n\ge0$)
$$-\sum_{v\in S}\sum_{k=1}^q\ln\|s_{kv}(P_n)\|_v\ge(q+2q\varepsilon)[K':\mathbb{Q}]{\rm h}_{\hat{L}}(P_n)-O(1)\quad.$$
D'où une contradiction avec $(1)$. $\square$\\

\subsection{Cas analytique}

Soit $X$ une variété complexe projective de dimension $d\ge1$.\\

{\bf Proposition 3.4:} {\it Soit $L\in{\rm Pic}(X)$ libre et gros.
Notons $q=h^0(X;L)$. On munit $L$ d'une métrique $\|\ \|$. Soient
$s_1;\cdots;s_N$ des sections non nulles engendrant $\Gamma(X;L)$. Soit
$\varepsilon>0$. Il existe alors un fermé $Z\neq X$ tel que pour toute courbe
entière $f$ sur $X$ d'image non contenue dans $Z(\mathbb{C})$, l'ensemble des
réels $r\ge0$ vérifiant
$$\int_0^{2\pi}\max_{J\in{\cal L}}\sum_{j\in J}\ln\|s_j(f(re^{i\theta}))\|^{-1}\frac{{\rm d}\theta}{2\pi}\ge(q+q\varepsilon){\rm T}_{\hat{L};f}(r)\qquad(1')$$
est de mesure de Lebesgue finie, où ${\cal L}$ désigne l'ensemble des parties
$J$ de $\{1;\cdots;n\}$ telles que $(s_j)_{j\in J}$ soit une base de
$\Gamma(X;L)$.}\\

{\it Démonstration:} En posant $V=\Gamma(X;L)$, on a un morphisme
$\Phi_L:X\rightarrow\mathbb{P}(V)$ génériquement fini. Il existe donc un fermé
$Z_1\neq X$ tel que $\Phi_{L|X-Z_1}$ soit à fibres finies. On applique alors la
version de Vojta ({\it cf} théorème 2 de \cite{Voj4}) du théorème de Cartan:

Il existe une réunion finie $H$ d'hyperplans de
$\mathbb{P}(V)\simeq\mathbb{P}^{q-1}_\mathbb{C}$ telle que pour toute courbe
entière d'image non contenue dans $Z_1\cup\Phi_L^{-1}(H)$, l'ensemble des réels
$r\ge0$ vérifiant $(1')$ est de mesure de Lebesgue finie. $\square$\\

{\bf Théorème 3.5:} {\it Soient $D_1;\cdots;D_r$ des diviseurs effectifs non
nuls sur $X$ qui se coupent proprement deux à deux. Posons
$Y=X-D_1\cup\cdots\cup D_r$. Soit $m\ge1$ un entier. On suppose que le diviseur
$\D L=m\sum_{i=1}^rD_i$ est libre et gros sur $X$ et que
$\nu(L;D_1;\cdots;D_r)>m$. Alors $Y$ est Brody quasi-hyperbolique.}\\

{\it Démonstration:} On procède en deux étapes: dans la première, on construit
un fermé $Z\neq X$ candidat à contenir toutes les courbes entières; dans la
seconde, on prouve que $Y$ est Brody quasi-hyperbolique.\\

{\it \'Etape 1:} On reprend la démonstration du théorème 3.3, jusqu'à la
construction des bases ${\cal B}_{I;\underline{a}}$. On munit chaque faisceau
${\cal O}_X(D_i)$ d'une métrique $\|\ \|$. Appliquons le théorème de Cartan et
Vojta (proposition 3.4) avec $\D\{s_1;\cdots;s_N\}={\cal B}_0\cup\bigcup_{I;\underline{a}}{\cal B}_{I;\underline{a}}$:

Il existe un fermé $Z\neq X$ tel que pour toute courbe entière $f$ sur $X$
d'image non contenue dans $Z(\mathbb{C})$, l'ensemble des réels $r\ge0$
vérifiant l'inégalité $(1')$ est de mesure de Lebesgue finie.\\

{\it \'Etape 2:} Soit $f$ une courbe entière sur $Y$. Raisonnons par l'absurde
en supposant $f(\mathbb{C})\not\subset Z(\mathbb{C})$.

Par compacité de $X(\mathbb{C})$, il existe un réel $M>0$ tel que pour tout
$y\in X(\mathbb{C})$, l'ensemble d'indices
$I_y=\{i\in\{1;\cdots;r\}\ |\ -\ln\|1_{D_i}(y)\|\ge M\}$ appartienne à
${\cal P}\cup\{\emptyset\}$ (il suffit d'extraire du recouvrement ouvert
$\Bigl(\Bigl\{y\in X(\mathbb{C})\ \Big|\ \exists I\in{\cal P}\cup\{\emptyset\}\ \forall i\notin I\ -\ln\|1_{D_i}(y)\|<M\Bigr\}\Bigr)_{M>0}$ un recouvrement
fini).\\

\underline{Fait}: Soit $y\in Y(\mathbb{C})$. Il existe une base
$(s_{1y};\cdots;s_{qy})$ de $\Gamma(X;L)$ contenue dans $\{s_1;\cdots;s_N\}$
telle que l'on ait la minoration suivante:
$$-\sum_{k=1}^q\ln\|s_{ky}(y)\|\ge-(q+3q\varepsilon)\ln\|1_L(y)\|-O(1)\quad,\qquad(2')$$
où le $O(1)$ est indépendant de $y$.\\

Prouvons ce fait. Si $I_y$ est vide, on prend
$\{s_{1y};\cdots;s_{qy}\}={\cal B}_0$ et on obtient la minoration $(2')$ en
remarquant que $-\ln\|1_L(y)\|<Mr$.

On suppose maintenant $I_y$ non vide. On a donc $I_y\in{\cal P}$. Pour tout
$i\in I_y$, on pose
$\D t_{yi}=\frac{\ln\|1_{D_i}(y)\|}{\D\sum_{j\in I_y}\ln\|1_{D_j}(y)\|}$.
Remarquons que l'on a $\D\sum_{i\in I_y} t_{yi}=1$. Choisissons un
$\underline{a}_y=(a_{yi})_i\in\Delta_{I_y}$ tel que $|bt_{yi}-a_{yi}|\le1$
pour tout $i\in I_y$. On prend alors
$\{s_{1y};\cdots;s_{qy}\}={\cal B}_{I_y;\underline{a}_y}$. Vérifions que ce
choix convient.\\

Soit $s\in\Gamma(X;L)-\{0\}$. Puisque les diviseurs $D_i$ se coupent proprement
deux à deux, le diviseur $\D{\rm div}(s)-\sum_{i\in I_y}\mu_i(s)D_i$ est
effectif. Ceci implique
$$-\ln\|s(y)\|\ge-\sum_{i\in I_y}\mu_i(s)\ln\|1_{D_i}(y)\|-O(1)\quad.$$

En remarquant que $\D t_{yi}\ge\frac{a_{yi}}{b}-\frac{m\varepsilon}{cr}$ et que
$\mu_i(s)\le c$ pour tout $i\in I_y$, on a, par définition des $t_{yi}$:
$$-\mu_i(s)\ln\|1_{D_i}(y)\|\ge-\Bigl(\frac{a_{yi}}{b}\mu_i(s)-\frac{m\varepsilon}{r}\Bigr)\sum_{j\in I_y}\ln\|1_{D_j}(y)\|$$
pour tout $i\in I_y$.

On en déduit l'inégalité
$$-\ln\|s(y)\|\ge-\Bigl(\frac{1}{b}\sum_{i\in I_y}a_{yi}\mu_i(s)-m\varepsilon\Bigr)\sum_{j\in I_y}\ln\|1_{D_j}(y)\|-O(1)\quad.$$

On écrit cette inégalité pour $s=s_{ky}$, puis on somme sur $k$. En observant
que pour $i\in I_y$, on a
$$\sum_{k=1}^q\sum_{i\in I_y}a_{yi}\mu_i(s_{ky})\ge\nu(L;D_1;\cdots;D_r)qb=(1+4\varepsilon)qbm$$
comme dans la démonstration du théorème 3.3, on trouve alors
$$-\sum_{k=1}^q\ln\|s_{ky}(y)\|\ge-(q+3q\varepsilon)m\sum_{j\in I_y}\ln\|1_{D_j}(y)\|-O(1)\quad.$$

Le fait énoncé $(2')$ s'en déduit en remarquant que $-\ln\|1_{D_j}(y)\|<M$
pour tout $j\notin I_y$.\\

Maintenant $f$ est une courbe entière sur $Y$, donc pour tout $r\ge0$, on a
$${\rm T}_{\hat{L};f}(r)=-\int_0^{2\pi}\ln\|1_L(f(re^{i\theta}))\|\frac{{\rm d}\theta}{2\pi}+O(1)\quad.$$

En utilisant la minoration $(2')$, on obtient (pour tout $r\ge0$)
$$\int_0^{2\pi}\max_{J\in{\cal L}}\sum_{k\in J}\ln\|s_k(f(re^{i\theta}))\|^{-1}\frac{{\rm d}\theta}{2\pi}\ge(q+3q\varepsilon){\rm T}_{\hat{L};f}(r)-O(1)\quad.$$
D'où une contradiction avec $(1')$ (puisque ${\rm T}_{\hat{L};f}(r)$ tend vers
$+\infty$ lorsque $r$ tend vers $+\infty$). $\square$\\

\section{Démonstration du théorème 2.1}

Soient $K$ un corps de caractéristique nulle et $X$ une variété projective sur
$K$ de dimension $d\ge2$. Pour tous diviseurs $L_1;\cdots;L_d$ sur $X$, on
désigne par $\bigl<L_1\cdots L_d\bigr>$ leur nombre d'intersection. Lorsque
$L$ est un diviseur sur $X$ tel que $q=h^0(X;L)\ge1$ et $E$ un diviseur
effectif non nul sur $X$, on pose
$\D\alpha(L;E)=\frac{1}{q}\sum_{k\ge1}h^0(X;L-kE)$.\\

{\bf Proposition 4.1:} {\it Soit $L$ un diviseur sur $X$ tel que
$q=h^0(X;L)\ge1$. Soient $D_1;\cdots;D_r$ des diviseurs effectifs non nuls sur
$X$ qui se coupent proprement deux à deux; supposons que toute intersection de
$\delta+1$ quelconques d'entre eux est vide (avec $2\le\delta\le r$). On a
alors
$$\nu(L;D_1;\cdots;D_r)\ge\frac{2}{\delta}\inf_i\alpha(L;D_i)\quad.$$}\\

{\it Démonstration:} On utilise les notations de la section 2.1. Lorsque $x$
est un réel, on désigne par $\lceil x\rceil$ le plus petit entier $\ge x$.
Soient $I\in{\cal P}$ et $\underline{a}\in\mathbb{N}^I-\{0\}$. Quitte à réduire
$I$, on peut supposer que $a_i\ge1$ pour tout $i\in I$. Observons que
$\# I\le\delta$.

Si $I$ est un singleton $\{i\}$, alors $\D V_{I;\underline{a};k}=\Gamma\Bigl(X;L-\Bigl\lceil\frac{k}{a_i}\Bigr\rceil D_i\Bigr)$,
donc on a bien
$$\sum_{k\ge1}V_{I;\underline{a};k}=a_i\sum_{k\ge1}h^0(X;L-kD_i)\ge a_iq\frac{2}{\delta}\inf_j\alpha(L;D_j)\quad.$$

On suppose maintenant $\# I\ge2$. On choisit deux indices $j<l$ dans $I$ tels
que $a_l\ge a_j\ge a_i$ pour tout $i\in I-\{j;l\}$. Pour
$(b_1;b_2)\in\mathbb{N}^2$, on pose $W(b_1;b_2)=\Gamma(X;L-b_1D_j-b_2D_l)$.

Soit $k$ un entier $\ge1$. L'espace vectoriel $V_{I;\underline{a};k}$ contient
alors le sous-espace
$$V'_k=W\Bigl(0;\Bigl\lceil\frac{k}{a_l}\Bigr\rceil\Bigr)+\sum_{b=0}^{\lceil k/a_l\rceil-1}W\Bigl(\Bigl\lceil\frac{k-a_lb}{a_j}\Bigr\rceil;b\Bigr)\quad.$$

Puisque les diviseurs $D_j$ et $D_l$ se coupent proprement, on a l'égalité
suivante pour tout $b'\in\{0;\cdots;\lceil k/a_l\rceil-1\}$:
$$W\Bigl(\Bigl\lceil\frac{k-a_lb'}{a_j}\Bigr\rceil;b'\Bigr)\bigcap\Bigl[W\Bigl(0;\Bigl\lceil\frac{k}{a_l}\Bigr\rceil\Bigr)+\sum_{b=b'+1}^{\lceil k/a_l\rceil-1}W\Bigl(\Bigl\lceil\frac{k-a_lb}{a_j}\Bigr\rceil;b\Bigr)\Bigr]=W\Bigl(\Bigl\lceil\frac{k-a_lb'}{a_j}\Bigr\rceil;b'+1\Bigr)\quad.$$

En utilisant $\lceil k/a_l\rceil$ fois la formule
$\dim(W_1+W_2)=\dim W_1+\dim W_2-\dim W_1\cap W_2$, on obtient que la dimension
de $V'_k$ vaut
$$\dim W\Bigl(0;\Bigl\lceil\frac{k}{a_l}\Bigr\rceil\Bigr)+\sum_{b=0}^{\lceil k/a_l\rceil-1}\Bigl[\dim W\Bigl(\Bigl\lceil\frac{k-a_lb}{a_j}\Bigr\rceil;b\Bigr)-\dim W\Bigl(\Bigl\lceil\frac{k-a_lb}{a_j}\Bigr\rceil;b+1\Bigr)\Bigr]\ .$$

Maintenant, on somme sur $k$ l'égalité précédente. Après simplifications, on
trouve
$$\sum_{k\ge1}\dim V'_k=a_l\sum_{k\ge1}h^0(X;L-kD_l)+a_j\sum_{k\ge1}h^0(X;L-kD_j)\quad.$$

On en conclut la minoration
$$\sum_{k\ge1}\dim V_{I;\underline{a};k}\ge\sum_{k\ge1}\dim V'_k\ge(a_j+a_l)q\inf_i\alpha(L;D_i)\ge\Bigl(\sum_{i\in I}a_i\Bigr)q\frac{2}{\delta}\inf_j\alpha(L;D_j)\quad.$$

D'où le résultat. $\square$\\

On aura besoin dans la suite d'une variante des ``inégalités de Morse
holomorphes'' ({\it cf} \cite{Dema} \S 12 et \cite{Ange}):\\

{\bf Lemme 4.2:} {\it Soient $E$ un diviseur libre et gros sur $X$ et $L$ un
diviseur sur $X$ tel que $L-E$ soit nef. Soit $\beta$ un réel $>0$. Pour tout
couple d'entiers $(n;k)$ vérifiant $1\le k\le\beta n$, on a alors la minoration
$$h^0(X;nL-kE)\ge\frac{\bigl<L^d\bigr>}{d!}n^d-\frac{\bigl<L^{d-1}E\bigr>}{(d-1)!}n^{d-1}k+\frac{d-1}{d!}\bigl<L^{d-2}E^2\bigr>n^{d-2}\min(k^2;n^2)-O(n^{d-1})\ ,$$
où le $O$ ne dépend pas de $(n;k)$.}\\

{\it Démonstration:} On a deux cas.\\

\underline{Cas $k\le n$}: La formule de Hirzebruch-Riemann-Roch donne que
$\chi(X;nL-kE)$ est une fonction polynomiale en $(n;k)$ dont on peut expliciter
la composante homogène dominante:

Pour tout $(n;k)$ tel que $1\le k\le n$, on a
$\D\chi(X;nL-kE)=\frac{1}{d!}\bigl<(nL-kE)^d\bigr>+O(n^{d-1})$.\\

Par ailleurs, d'après le théorème 1.4.40 de \cite{Laza} p. 69 (ou plutôt
d'après sa démonstration), on a $h^i(X;nL-kE)=O(n^{d-i})$ pour tout $i\ge1$,
puisque $L$ et $L-E$ sont nefs. On a en particulier
$\D h^0(X;nL-kE)=\frac{1}{d!}\bigl<(nL-kE)^d\bigr>+O(n^{d-1})$.\\

Or un calcul montre (par multilinéarité) la formule
$$\bigl<(nL-kE)^d\bigr>=\bigl<L^d\bigr>n^d-d\bigl<L^{d-1}E\bigr>n^{d-1}k+\sum_{i=2}^d(i-1)\bigl<L^{i-2}(nL-kE)^{d-i}E^2\bigr>n^{i-2}k^2\quad.$$

L'inégalité de l'énoncé s'en déduit facilement: les diviseurs $L$, $nL-kE$ et
$E$ sont nefs, donc on a $\bigl<L^{i-2}(nL-kE)^{d-i}E^2\bigr>\ge0$ pour tout
$i\in\{2;\cdots;d-1\}$.\\

\underline{Cas $k>n$}: D'après le théorème de Bertini ({\it cf} corollaire 6.11
de \cite{Joua} p. 89), il existe $s\in\Gamma(X;E)-\{0\}$ tel que
$Z={\rm div}(s)$ soit géométriquement intègre sur $K$.\\

Soit $i$ un entier tel que $n\le i\le\beta n$. On a la suite exacte de
${\cal O}_X$-modules suivante:
$$0\rightarrow{\cal O}_X(nL-(i+1)E)\rightarrow{\cal O}_X(nL-iE)\rightarrow{\cal O}_X(nL-iE)_{|Z}\rightarrow0\quad.$$

On en déduit une suite exacte en cohomologie qui fournit l'inégalité
$$h^0(X;nL-(i+1)E)\ge h^0(X;nL-iE)-h^0(Z;(nL-iE)_{|Z})\quad.$$

En utilisant la majoration
$$h^0(Z;(nL-iE)_{|Z})\le h^0(Z;nL_{|Z})=\frac{\bigl<L^{d-1}E\bigr>}{(d-1)!}n^{d-1}+O(n^{d-2})$$
(obtenue par Hirzebruch-Riemann-Roch), on trouve
$$\begin{array}{rcl}
\D h^0(X;nL-kE)&\ge&\D h^0(X;nL-nE)-\sum_{i=n}^{k-1}h^0(Z;(nL-iE)_{|Z})\\
&\ge&\D\frac{\bigl<L^d\bigr>}{d!}n^d-\frac{\bigl<L^{d-1}E\bigr>}{(d-1)!}n^{d-1}k+\frac{d-1}{d!}\bigl<L^{d-2}E^2\bigr>n^d-O(n^{d-1})\\
\end{array}$$
(la minoration de $h^0(X;nL-nE)$ est donnée par le premier cas). D'où le
résultat. $\square$\\

{\bf Remarque:} La démonstration fournit en fait une minoration de
$h^0(X;nL-kE)-h^1(X;nL-kE)$.\\

On note ici $g:\mathbb{R}_+\rightarrow\mathbb{R}_+$ l'application continue
définie par $\D g(\beta)=\frac{\beta^3}{3}$ si $\beta\le1$ et
$\D g(\beta)=\beta-\frac{2}{3}$ si $\beta\ge1$.\\

{\bf Corollaire 4.3:} {\it Soient $E$ un diviseur effectif libre et gros sur
$X$ et $L$ un diviseur sur $X$ tel que $L-E$ soit nef. On a alors
$$\liminf_{n\rightarrow+\infty}\frac{1}{n}\alpha(nL;E)\ge\frac{\bigl<L^d\bigr>}{2d\bigl<L^{d-1}E\bigr>}+(d-1)\frac{\bigl<L^{d-2}E^2\bigr>}{\bigl<L^d\bigr>}g\Bigl(\frac{\bigl<L^d\bigr>}{d\bigl<L^{d-1}E\bigr>}\Bigr)\quad.$$}\\

{\it Démonstration:} On pose
$\D\beta=\frac{\bigl<L^d\bigr>}{d\bigl<L^{d-1}E\bigr>}$ et
$M=(d-1)\bigl<L^{d-2}E^2\bigr>$. Grâce au lemme 4.2, on a les estimations
suivantes:
$$\begin{array}{rcl}
\D\sum_{k\ge1}h^0(X;nL-kE)&\ge&\D\sum_{k=1}^{\lfloor\beta n\rfloor}\Bigl(\frac{\bigl<L^d\bigr>}{d!}n^d-\frac{\bigl<L^{d-1}E\bigr>}{(d-1)!}n^{d-1}k+\frac{M}{d!}n^{d-2}\min(k^2;n^2)\Bigr)-O(n^d)\\
&=&\D\Bigl(\frac{\bigl<L^d\bigr>}{d!}\beta-\frac{\bigl<L^{d-1}E\bigr>}{(d-1)!}\frac{\beta^2}{2}+\frac{M}{d!}g(\beta)\Bigr)n^{d+1}-O(n^d)\quad.\\
\end{array}$$

D'où la minoration $\D\alpha(nL;E)\ge\Bigl(\frac{\beta}{2}+\frac{M}{\bigl<L^d\bigr>}g(\beta)\Bigr)n-O(1)$. $\square$\\

Montrons maintenant le résultat principal de cette section:\\

{\bf Théorème 4.4:} {\it Soient $D_1;\cdots;D_r$ des diviseurs effectifs
presque amples sur $X$. Il existe alors des entiers $m_1;\cdots;m_r$ tels qu'en
posant $\D L=\sum_{i=1}^rm_iD_i$, on ait
$$m_i\ge1\ \mbox{et}\ \liminf_{n\rightarrow+\infty}\frac{1}{n}\alpha(nL;m_iD_i)>\frac{r}{2d}\ \mbox{pour tout}\ i\in\{1;\cdots;r\}\ .$$}\\

{\it Démonstration:} On pose ici
$\Delta=\{(t_1;\cdots;t_r)\in\mathbb{R}_+^r\ |\ t_1+\cdots+t_r=1\}$.  Pour tout
$t=(t_1;\cdots;t_r)\in\Delta$, on désigne par $L_t$ le $\mathbb{R}$-diviseur
$\D L_t=\sum_{j=1}^rt_jD_j$ et on pose
$\D\phi(t)=\Bigl(\sum_{i=1}^r\frac{1}{\bigl<L_t^{d-1}D_i\bigr>}\Bigr)^{-1}$.

On note $f:\Delta\rightarrow\Delta$ l'application continue définie par
$\D f(t)=\Bigl(\frac{\phi(t)}{\bigl<L_t^{d-1}D_1\bigr>};\cdots;\frac{\phi(t)}{\bigl<L_t^{d-1}D_r\bigr>}\Bigr)$ pour tout $t\in\Delta$. D'après le théorème de
Brouwer, $f$ admet un point fixe $x=(x_1;\cdots;x_r)$. On a alors
$\phi(x)=\bigl<L_x^{d-1}D_i\bigr>x_i$ pour tout
$i\in\{1;\cdots;r\}$, donc $\phi(x)r=\bigl<L_x^d\bigr>$.\\

On en déduit l'inégalité
$$\frac{\bigl<L_x^d\bigr>}{2d\bigl<L_x^{d-1}D_i\bigr>x_i}+(d-1)\frac{\bigl<L_x^{d-2}D_i^2\bigr>x_i^2}{\bigl<L_x^d\bigr>}g\Bigl(\frac{\bigl<L_x^d\bigr>}{d\bigl<L_x^{d-1}D_i\bigr>x_i}\Bigr)>\frac{r}{2d}\ \mbox{pour tout}\ i\in\{1;\cdots;r\}\ .$$

On approche $x$ par un $y\in\mathbb{Q}_+^{*r}\cap\Delta$ de la forme
$\D y=\Bigl(\frac{m_1}{m};\cdots;\frac{m_r}{m}\Bigr)$ de telle sorte que
l'inégalité précédente soit encore valable avec $y$ au lieu de $x$, et on
conclut en appliquant le corollaire 4.3. $\square$\\

On en déduit le théorème 2.1 en appliquant la proposition 4.1.\\

\section{Géométrie bis}

\subsection{Préliminaires}

Soient $r$ et $m$ des entiers $\ge1$. On pose $\Delta=\{0;\cdots;m\}^r$. On
munit $\Delta$ de l'ordre lexicographique. Notons $\underline{m}=(m;\cdots;m)$
le plus grand élément de $\Delta$. Pour tout
$\underline{b}=(b_1;\cdots;b_r)\in\Delta$, on désigne par
$J_{\underline{b}}$ l'ensemble des $i\in\{1;\cdots;r\}$ tels que $b_i<m$.\\

Commençons par la variante suivante du lemme 2.2 de \cite{CZ1}:\\

{\bf Lemme 5.1:} {\it Soit $A$ un anneau local. Soit
$(\varphi_1;\cdots;\varphi_r)$ une suite régulière de $A$. Pour tout
$\underline{b}\in\Delta$, on a alors l'inclusion d'idéaux
$$(\varphi_1^{b_1}\cdots\varphi_r^{b_r}A)\cap\Bigl(\sum_{\underline{c}>\underline{b}}\varphi_1^{c_1}\cdots\varphi_r^{c_r}A\Bigr)\subset\sum_{j\in J_{\underline{b}}}\varphi_1^{b_1}\cdots\varphi_r^{b_r}\varphi_jA\quad.$$}\\

{\it Démonstration:} On raisonne par récurrence sur $r$. Si $r=1$, alors
l'inclusion est évidente. Supposons $r\ge2$ et le résultat au cran $r-1$.
Posons $\Delta'=\{0;\cdots;m\}^{r-1}$ et $\underline{b}'=(b_2;\cdots;b_r)$.
Soit $\D x\in(\varphi_1^{b_1}\cdots\varphi_r^{b_r}A)\cap\Bigl(\sum_{\underline{c}>\underline{b}}\varphi_1^{c_1}\cdots\varphi_r^{c_r}A\Bigr)$. On a deux cas.\\

\underline{Cas $b_1=m$}: L'élément $x$ s'écrit $\D x=\varphi_1^{m}y=\sum_{\underline{c}'>\underline{b}'}\varphi_1^{m}a_{\underline{c}'}$ avec un
$y\in\varphi_2^{b_2}\cdots\varphi_r^{b_r}A$ et des
$a_{\underline{c}'}\in\varphi_2^{c_2}\cdots\varphi_r^{c_r}A$. En simplifiant
par $\varphi_1^m$, on obtient que $y$ appartient à $\D\sum_{\underline{c}'>\underline{b}'}\varphi_2^{c_2}\cdots\varphi_r^{c_r}A$. Or
$(\varphi_2;\cdots;\varphi_r)$ est une suite régulière de $A$, donc $y$ est un
élément de $\D\sum_{j\in J_{\underline{b}}}\varphi_2^{b_2}\cdots\varphi_r^{b_r}\varphi_jA$ par hypothèse de récurrence.\\

\underline{Cas $b_1<m$}: L'élément $x$ s'écrit
$\D x=\varphi_1^{b_1}y=\varphi_1^{b_1+1}z+\sum_{\underline{c}'>\underline{b}'}\varphi_1^{b_1}a_{\underline{c}'}$ avec un
$y\in\varphi_2^{b_2}\cdots\varphi_r^{b_r}A$, un $z\in A$ et des
$a_{\underline{c}'}\in\varphi_2^{c_2}\cdots\varphi_r^{c_r}A$. On écrit
$y=\varphi_2^{b_2}\cdots\varphi_r^{b_r}w$ avec $w\in A$. On simplifie par
$\varphi_1^{b_1}$ puis on réduit modulo $\varphi_1$; on trouve ainsi dans
$A'=A/\varphi_1A$ l'égalité
$\D\bar{y}=\sum_{\underline{c}'>\underline{b}'}\bar{a_{\underline{c}'}}$.

On en déduit que $\bar{y}$ appartient à
$\D(\bar{\varphi_2}^{b_2}\cdots\bar{\varphi_r}^{b_r}A')\cap\Bigl(\sum_{\underline{c}'>\underline{b}'}\bar{\varphi_2}^{c_2}\cdots\bar{\varphi_r}^{c_r}A'\Bigr)$.
Or $(\bar{\varphi_2};\cdots;\bar{\varphi_r})$ est une suite régulière de $A'$,
donc $\bar{y}$ est un élément de $\D\sum_{j\in J_{\underline{b}}-\{1\}}\bar{\varphi_2}^{b_2}\cdots\bar{\varphi_r}^{b_r}\bar{\varphi_j}A'$  par hypothèse de
récurrence. En simplifiant par
$\bar{\varphi_2}^{b_2}\cdots\bar{\varphi_r}^{b_r}$, on obtient que $\bar{w}$
est dans $\D\sum_{j\in J_{\underline{b}}-\{1\}}\bar{\varphi_j}A'$. On en
conclut que $w$ appartient à $\D\sum_{j\in J_{\underline{b}}}\varphi_jA$.\\

D'où le résultat. $\square$\\

Soient $K$ un corps de caractéristique nulle et $X$ une variété projective sur
$K$ de dimension $d\ge1$.\\

{\bf Définition:} Un ${\cal O}_X$-module cohérent ${\cal C}$ sur $X$ est dit
{\bf acyclique} lorsque $h^i(X;{\cal C})=0$ pour tout $i\ge1$.\\

Soit $L$ un diviseur sur $X$ tel que $q=h^0(X;L)\ge1$. Soient $D_1;\cdots;D_r$
des diviseurs effectifs non nuls sur $X$ qui se coupent proprement.

Pour tout $\underline{b}\in\Delta$, on pose $\D{\cal L}_{\underline{b}}={\cal O}_X\Bigl(L-\sum_{i=1}^rb_iD_i\Bigr)$. Pour
$\underline{b}\in\Delta$, on définit le sous-module ${\cal C}_{\underline{b}}$
de ${\cal L}_{\underline{b}}$ par
$${\cal C}_{\underline{b}}=\sum_{j\in J_{\underline{b}}}{\cal O}_X\Bigl(L-D_j-\sum_{i=1}^rb_iD_i\Bigr)\quad.$$

Soit $(a_1;\cdots;a_r)\in\mathbb{N}^r$. Pour $k\in\mathbb{N}^*$, on note $V_k$
le sous-espace de $\Gamma(X;L)$ défini par
$$V_k=\sum_{\underline{b}}\Gamma(X;{\cal L}_{\underline{b}})\ \mbox{où la somme porte sur les}\ \underline{b}\in\mathbb{N}^r\ \mbox{tels que}\ \sum_{i=1}^ra_ib_i\ge k\ .$$

{\bf Lemme 5.2:} {\it Avec ces notations, on a la minoration suivante:
$$\sum_{k\ge1}\dim V_k\ge\sum_{i=1}^ra_i\sum_{\underline{b}\in\Delta}\Bigl[h^0(X;{\cal L}_{\underline{b}})-h^0(X;{\cal C}_{\underline{b}})\Bigr]b_i\quad.$$}\\

{\it Démonstration:} Soit $k$ un entier tel que $\D1\le k\le\sum_{i=1}^ra_im$.
Notons ${\cal D}_k$ l'ensemble des $\underline{b}\in\Delta$ tels que
$\D\sum_{i=1}^ra_ib_i\ge k$. L'espace vectoriel $V_k$ contient alors le
sous-espace $\D V'_k=\sum_{\underline{b}\in{\cal D}_k}\Gamma(X;{\cal L}_{\underline{b}})$.\\

Soit $\underline{b}\in{\cal D}_k-\{\underline{m}\}$. Le lemme 5.1 fournit
l'inclusion de ${\cal O}_X$-modules $\D{\cal L}_{\underline{b}}\cap\sum_{\underline{c}>\underline{b}}{\cal L}_{\underline{c}}\subset{\cal C}_{\underline{b}}$,
puisque les diviseurs $D_1;\cdots;D_r$ se coupent proprement. On a en
particulier l'inclusion d'espaces vectoriels
$$\Gamma(X;{\cal L}_{\underline{b}})\cap\sum_{\underline{c}>\underline{b}}\Gamma(X;{\cal L}_{\underline{c}})\subset\Gamma(X;{\cal C}_{\underline{b}})\quad.$$
En utilisant $\#{\cal D}_k-1$ fois la formule
$\dim(W_1+W_2)=\dim W_1+\dim W_2-\dim W_1\cap W_2$, on trouve l'inégalité
$$\dim V'_k\ge h^0(X;{\cal L}_{\underline{m}})+\sum_{\underline{b}\in{\cal D}_k-\{\underline{m}\}}\Bigr[h^0(X;{\cal L}_{\underline{b}})-h^0(X;{\cal C}_{\underline{b}})\Bigr]\quad.$$

On obtient le résultat en sommant sur $k$ cette inégalité. $\square$\\

{\bf Proposition 5.3:} {\it On suppose de plus que ${\cal L}_{\underline{b}}$
est acyclique pour tout $\underline{b}\in\Delta$. On a alors
$$\sum_{k\ge1}\dim V_k\ge\sum_{i=1}^ra_i\sum_{k=1}^mh^0(X;L-kD_i)\quad.$$
On a en particulier
$\D\nu(L;D_1;\cdots;D_r)\ge\frac{1}{q}\inf_i\sum_{k=1}^mh^0(X;L-kD_i)$.}\\

{\it Démonstration:} Soit $\underline{b}\in\Delta-\{\underline{m}\}$. Pour
toute partie $I$ de $J_{\underline{b}}$, posons ici
$${\cal L}_{\underline{b};I}={\cal O}_X\Bigl(L-\sum_{j\in J_{\underline{b}}}D_j-\sum_{i=1}^rb_iD_i\Bigr)\quad.$$
On pose aussi $p=\#J_{\underline{b}}$ et $\D{\cal E}_{\underline{b}}=\bigoplus_{j\in J_{\underline{b}}}{\cal L}_{\underline{b};\{j\}}$.

Les diviseurs $(D_j)_{j\in J_{\underline{b}}}$ se coupent proprement, donc on a
la suite exacte de Koszul suivante ({\it cf} \cite{Fult} p. 431):
$$0\rightarrow\Lambda^p{\cal E}_{\underline{b}}\rightarrow\cdots\rightarrow\Lambda^1{\cal E}_{\underline{b}}\rightarrow{\cal C}_{\underline{b}}\rightarrow0\quad.$$

On remarque que $\D\Lambda^j{\cal E}_{\underline{b}}=\bigoplus_{\#I=j}{\cal L}_{\underline{b};I}$ (qui est en particulier acyclique) pour tout
$j\in\{1;\cdots;p\}$. La suite de Koszul précédente induit donc par acyclicité
une suite exacte en image directe:
$$0\rightarrow\Gamma(X;\Lambda^p{\cal E}_{\underline{b}})\rightarrow\cdots\rightarrow\Gamma(X;\Lambda^1{\cal E}_{\underline{b}})\rightarrow\Gamma(X;{\cal C}_{\underline{b}})\rightarrow0\quad.$$

On en déduit la relation
$$h^0(X;{\cal L}_{\underline{b}})-h^0(X;{\cal C}_{\underline{b}})=\sum_{I\subset J_{\underline{b}}}(-1)^{\#I}h^0(X;{\cal L}_{\underline{b};I})\quad.$$

Maintenant, on fixe $i\in\{1;\cdots;r\}$ et $c\in\{0;\cdots;m\}$, et on somme
sur l'ensemble $\Delta'_c$ des $\underline{b}\in\Delta$ tels que $b_i=c$. Un
réarrangement des termes permet de simplifier et montre que:
$$\sum_{\underline{b}\in\Delta'_c}\sum_{I\subset J_{\underline{b}}}(-1)^{\#I}h^0(X;{\cal L}_{\underline{b};I})=\begin{cases}
h^0(X;L-cD_i)-h^0(X;L-(c+1)D_i)&\mbox{si}\ c<m;\\
h^0(X;L-mD_i)&\mbox{si}\ c=m.\\
\end{cases}$$
(En effet, si $p'=\#\{j\neq i\ |\ b_j\ge1\}\ge1$, alors le terme
$h^0(X;{\cal L}_{\underline{b}})$ apparaît $2^{p'-1}$ fois avec le signe plus
et $2^{p'-1}$ fois avec le signe moins; de même avec le terme
$h^0(X;{\cal L}_{\underline{b};\{i\}})$ dans le cas $c<m$).\\

On en déduit l'égalité
$$\begin{array}{rcl}
\D\sum_{\underline{b}\in\Delta}\Bigl[h^0({\cal L}_{\underline{b}})-h^0({\cal C}_{\underline{b}})\Bigr]b_i&=&\D h^0(L-mD_i)m+\sum_{c=0}^{m-1}\Bigl[h^0(L-cD_i)-h^0(L-(c+1)D_i)\Bigr]c\\
&=&\D\sum_{k=1}^mh^0(X;L-kD_i)\quad.\\
\end{array}$$

On conclut en appliquant le lemme 5.2. $\square$\\

\subsection{Démonstration du théorème 2.2}

Soient $K$ un corps de caractéristique nulle et $X$ une variété projective sur
$K$ de dimension $d\ge1$. Soient $D_1;\cdots;D_r$ des diviseurs effectifs non
nuls sur $X$ qui se coupent proprement (avec $r\ge1$). Notons ${\cal P}$
l'ensemble des parties $I$ non vides de $\{1;\cdots;r\}$ telles que
$\D\bigcap_{i\in I}D_i$ soit non vide.\\

{\bf Théorème 5.4:} {\it Soit $L$ un diviseur ample sur $X$. On suppose que
$D_i$ est nef pour tout $i\in\{1;\cdots;r\}$. Soit $\theta>1$ un réel tel que
le $\mathbb{R}$-diviseur $\D L-\theta\sum_{i\in I}D_i$ soit nef pour tout
$I\in{\cal P}$. On a alors l'inégalité
$$\liminf_{n\rightarrow+\infty}\frac{1}{n}\nu(nL;D_1;\cdots;D_r)\ge\frac{\theta}{(d+1)\bigl<L^d\bigr>}\inf_i\sum_{j=0}^d\bigl<L^{d-j}(L-\theta D_i)^j\bigr>\quad.$$}\\

{\it Démonstration:} D'après le théorème d'annulation de Fujita ({\it cf}
théorème 1.4.35 de \cite{Laza} p. 66), il existe $n_0\ge1$ tel que $n_0L+N$
soit acyclique pour tout $N\in{\rm Pic}(X)$ nef.

On pose $n'=\lfloor(n-n_0)\theta\rfloor$ pour tout $n>n_0$. En appliquant la
proposition 5.3 (avec $m=n'$), on obtient (pour tout $n>n_0$)
$$\nu(nL;D_1;\cdots;D_r)\ge\frac{1}{h^0(X;nL)}\inf_i\sum_{k=1}^{n'}h^0(X;nL-kD_i)\quad.$$

Soit $i\in\{1;\cdots;r\}$. Grâce à la formule de Hirzebruch-Riemann-Roch, on a
les estimations suivantes:
$$\begin{array}{rcl}
\D\sum_{k=1}^{n'}h^0(X;nL-kD_i)&=&\D\frac{1}{d!}\sum_{k=1}^{n'}\Bigl[\bigl<(nL-kD_i)^d\bigr>+O(n^{d-1})\Bigr]\\
&=&\D\frac{1}{d!}\sum_{j=0}^d\sum_{k=1}^{n'}{\rm C}_d^j\bigl<L^{d-j}D_i^j\bigr>n^{d-j}(-k)^j+O(n^d)\\
&=&\D\frac{1}{d!}\sum_{j=0}^d{\rm C}_d^j\bigl<L^{d-j}D_i^j\bigr>\frac{(-1)^j}{j+1}\theta^{j+1}n^{d+1}+O(n^d)\quad.\\
\end{array}$$

Or un calcul montre la formule
$$\sum_{j=0}^d{\rm C}_d^j\bigl<L^{d-j}D_i^j\bigr>\frac{(-1)^j}{j+1}\theta^{j+1}=\frac{\theta}{d+1}\sum_{j=0}^d\bigl<L^{d-j}(L-\theta D_i)^j\bigr>\quad.$$

D'où l'inégalité de l'énoncé. $\square$\\

{\bf Corollaire 5.5:} {\it Soit $L$ un diviseur ample sur $X$. On suppose que
$D_i$ est nef pour tout $i\in\{1;\cdots;r\}$. Soit $\theta>1$ un réel tel que
le $\mathbb{R}$-diviseur $\D L-d\theta D_i$ soit nef pour tout
$i\in\{1;\cdots;r\}$. On a alors
$$\liminf_{n\rightarrow+\infty}\frac{1}{n}\nu(nL;D_1;\cdots;D_r)\ge\lambda_d\theta\quad.$$}\\

{\it Démonstration:} Pour tout $I\in{\cal P}$, le $\mathbb{R}$-diviseur
$\D L-\theta\sum_{i\in I}D_i$ est nef puisque $\#I\le d$. Soit
$i\in\{1;\cdots;r\}$. Les $\mathbb{R}$-diviseurs
$\D L-\theta D_i-\Bigl(1-\frac{1}{d}\Bigr)L$ et $L$ sont nefs, donc on a
$$\bigl<L^{d-j}(L-\theta D_i\bigr)^j\bigr>\ge\Bigl(1-\frac{1}{d}\Bigr)^j\bigl<L^d\bigr>\ \mbox{pour tout}\ j\in\{1;\cdots;d\}\ .$$

En appliquant le théorème 5.4, on trouve ainsi
$$\liminf_{n\rightarrow+\infty}\frac{1}{n}\nu(nL;D_1;\cdots;D_r)\ge\frac{\theta}{d+1}\sum_{j=0}^d\Bigl(1-\frac{1}{d}\Bigr)^j=\lambda_d\theta\quad.$$

D'où le résultat. $\square$\\

\ \\

{\small Pascal Autissier. I.R.M.A.R., Université de Rennes I, campus de
Beaulieu, 35042 Rennes cedex, France.

pascal.autissier@univ-rennes1.fr}

\end{document}